# Distributed Maximum Likelihood using Dynamic Average Consensus Algorithm

J. George

*Abstract*—This paper presents the formulation and analysis of a novel distributed maximum likelihood algorithm that utilizes a first-order optimization scheme. The proposed approach utilizes a static average consensus algorithm to reach agreement on the initial condition to the iterative optimization scheme and a dynamic average consensus algorithm to reach agreement on the gradient direction. The current distributed algorithm is guaranteed to exponentially recover the performance of the centralized algorithm. Though the current formulation focuses on maximum likelihood algorithm built on first-order methods, it can be easily extended to higher order methods. Numerical simulations validate the theoretical contributions of the paper.

*Index Terms*—Distributed maximum likelihood, sensor network, dynamic average consensus, first-order methods

## I. INTRODUCTION

The advent of sensor network for a broad range of surveillance and reconnaissance applications has highlighted the utility of scalable algorithms which can be implemented in a distributed fashion. Several distributed maximum likelihood algorithms have been proposed to tackle the parameter estimation problem in sensor network [1]–[4]. These algorithms typically utilize an existing distributed optimization scheme and are not focused on recovering the centralized estimator trajectory.

Maximum likelihood problems can be cast as optimization problems, and more specifically, in the context of sensor network, as distributed optimization problems. Early approaches to distributed optimization involve the Distributed Subgradient Methods (DSMs) [5], where the non-smooth function optimization is performed by means of subgradient-based descent or ascent approaches. Typical approaches to DSMs involve primal subgradient methods [6]–[8] that yield sublinear convergence rates. DSMs have been exploited for several practical purposes, e.g., to optimally allocate resources in Wireless Sensor Networks (WSNs) [9], to maximize the convergence speeds of gossip algorithms [10], and to manage optimality criteria defined in terms of ergodic limits [11]. While DSMs have the advantage of being easily distributed, have limited computational requirements, and are inherently asynchronous [12]–[14], they suffer from a low convergence rate since they require the update step size to decrease to zero as $k \to \infty$ ($k$ being the iteration step). Therefore, the rate of convergence is subexponential and methods like multi-step approaches [15] and Newton-like methods [16]–[18] have proposed to improve the convergence rate of DSMs. Examples of distributed maximum likelihood algorithms using DSM can be found in [19]–[21].

Among the distributed optimization methods, the most widely known algorithm is Alternating Direction Method of Multipliers (ADMM) [22], whose roots can be traced back to [23]. Recent advancements in ADMM for asynchronous and distributed implementations can be found in [24]–[27]. Though it is efficient in several practical scenarios [28], ADMM often requires the agents to reach consensus on the design variable at each iteration of the optimization step and does not offer the robustness necessary for sensor network. Examples of ADMM based distributed maximum likelihood methods include [3] and [4].

Recently, several alternative approaches to ADMM and DSM have appeared. For example, in [8], [29] the authors construct contraction mappings by means of cyclic projections of the estimate of the optimum onto the constraints. Other methods include the F-Lipschitz methods [30]; the distributed randomized Kaczmarz method [31]; Zero Gradient Sum (ZGS) algorithms [32]; exact first-order algorithm (EXTRA) [2]; and distributed dual subgradient methods [33].

While the above approaches focus on developing new distributed optimization algorithms, here we develop a framework for the distributed implementation of existing optimization algorithms for maximum likelihood estimation. Though we focus on firs-order mthods, the proposed approach can be easily utilized for the implementation of maximum likelihood estimation using higher-order

Jemin George is with the U.S. Army Research Laboratory, Adelphi, MD 20783, USA. `jemin.george.civ@mail.mil`

methods. Current approach utilizes a *dynamic average consensus*[1] (DAC) algorithm [34], [35] to reach agreement on the gradient direction in finite time. Unlike most of the existing distributed algorithms that only focus on asymptotically recovering the centralized steady-state solution, the proposed approach guarantees exponential convergence of the distributed estimator trajectories to the centralized estimator trajectories.

The rest of this paper is organized as follows. Mathematical preliminaries and detailed problem formulation are given in Sections II and III, respectively. Main results of the paper are given in Section IV. Section V provides the results obtained from numerical simulations. Conclusions and future work are discussed in Section VI

## II. PRELIMINARIES

### A. Notation

Let $\mathbb{R}^{n \times m}$ denote the set of $n \times m$ real matrices. An $n \times n$ identity matrix is denoted as $I_n$ and $\mathbf{1}_n$ denotes an $n$-dimensional vector of all ones. Let $\mathbb{R}_{\mathbf{1}}^n$ denote the set of all $n$-dimensional vectors of the form $\kappa \mathbf{1}_n$, where $\kappa \in \mathbb{R}$. The absolute value of a vector is given as $|\mathbf{x}| = \begin{bmatrix} |x_1| & \ldots & |x_n| \end{bmatrix}^T$. Let $\operatorname{sgn}\{\cdot\}$ denote the signum function, defined as

$$\operatorname{sgn}\{x\} \triangleq \begin{cases} +1, & \text{if } x > 0; \\ 0, & \text{if } x = 0; \\ -1, & \text{if } x < 0, \end{cases}$$

and $\forall \mathbf{x} \in \mathbb{R}^n$, $\operatorname{sgn}\{\mathbf{x}\} \triangleq \begin{bmatrix} \operatorname{sgn}\{x_1\} & \ldots & \operatorname{sgn}\{x_n\} \end{bmatrix}^T$. For $p \in [1, \infty]$, the $p$-norm of a vector $\mathbf{x}$ is denoted as $\|\mathbf{x}\|_p$. For matrices $A \in \mathbb{R}^{m \times n}$ and $B \in \mathbb{R}^{p \times q}$, $A \otimes B \in \mathbb{R}^{mp \times nq}$ denotes their Kronecker product.

### B. Network Model

For a *connected undirected* graph $\mathcal{G}(\mathcal{V}, \mathcal{E})$ of order $n$, $\mathcal{V} \triangleq \{v_1, \ldots, v_n\}$ represents the agents or nodes. The communication links between the agents are represented as $\mathcal{E} \triangleq \{e_1, \ldots, e_\ell\} \subseteq \mathcal{V} \times \mathcal{V}$. Here each undirected edge is considered as two distinct directed edges and the edges are labeled such that they are grouped into incoming links to nodes $v_1$ to $v_n$. Let $\mathcal{I}$ denote the index set $\{1, \ldots, n\}$ and $\forall i \in \mathcal{I}$; let $\mathcal{N}_i \triangleq \{v_j \in \mathcal{V} : (v_i, v_j) \in \mathcal{E}\}$ denote the set of neighbors of node $v_i$. Let $\mathcal{A} \triangleq [a_{ij}] \in \{0, 1\}^{n \times n}$ be the *adjacency matrix* with entries $a_{ij} = 1$ if $(v_i, v_j) \in \mathcal{E}$ and zero otherwise. Define $\Delta \triangleq \operatorname{diag}(\mathcal{A} \mathbf{1}_n)$ as the degree matrix associated with the graph and $\mathcal{L} \triangleq \Delta - \mathcal{A}$ as the graph *Laplacian*. The *incidence matrix* of the graph is defined as $\mathcal{B} = [b_{ij}] \in \{-1, 0, 1\}^{n \times \ell}$, where $b_{ij} = -1$ if edge $e_j$ leaves node $v_i$, $b_{ij} = 1$ if edge $e_j$ enters node $v_i$, and $b_{ij} = 0$ otherwise.

[1] The problem of reaching consensus on the average of a set of local time-varying signals in a distributed fashion is typically referred to as the *dynamic average consensus*.

For the connected undirected graph $\mathcal{G}(\mathcal{V}, \mathcal{E})$, $\mathcal{L} = \frac{1}{2}\mathcal{B}\mathcal{B}^T$ and $\mathcal{L}$ is a positive semi-definite matrix with one eigenvalue at 0 corresponding to the eigenvector $\mathbf{1}_n$. Furthermore, we have $M \triangleq \left( I_n - \frac{1}{n} \mathbf{1}_n \mathbf{1}_n^T \right) = \mathcal{L}(\mathcal{L})^+ = \mathcal{B}\mathcal{B}^T (\mathcal{B}\mathcal{B}^T)^+ = \mathcal{B} (\mathcal{B}^T \mathcal{B})^+ \mathcal{B}^T$, where $(\cdot)^+$ denotes the generalized inverse.

*Remark 1:* For all $\mathbf{x} \in \mathbb{R}^n$, such that $\mathbf{1}_n^T \mathbf{x} = 0$, we have $\mathbf{x}^T \mathcal{L} (\mathcal{L})^+ \mathbf{x} = \mathbf{x}^T \mathbf{x} > 0$.

## III. PROBLEM FORMULATION

Consider a network of sensors, represented as a connected, undirected network $\mathcal{G}(\mathcal{V}, \mathcal{E})$ of order $n$, where the nodes represents the sensors and the edges represent the communication links between the sensors. Sensors $v_i$ and $v_j$ are (one-hop) neighbors if $(v_i, v_j) \in \mathcal{E}$. We will assume that all sensor are synchronized with a common clock and each sensor can only communicate with its neighboring sensors. For all $i \in \mathcal{I}$, the individual sensor measurements are given as

$$\mathbf{z}_i = \mathbf{h}_i(\boldsymbol{\theta}) + \mathbf{w}_i, \quad (1)$$

where $\boldsymbol{\theta} \in \mathbb{R}^r$ is the latent variable to be estimated and $\mathbf{w}_i \in \mathbb{R}^m$ is the measurement noise associated with the $i^{\text{th}}$-sensor. Noise is assumed to be zero-mean, independent, Gaussian noise with known variance, i.e., $\mathbf{w}_i \sim \mathrm{N}(\mathbf{0}, R_i)$. The nonlinear mapping $\mathbf{h}_i(\cdot) : \mathbb{R}^r \mapsto \mathbb{R}^m$ is locally known to each sensors.

Under the current setup, the optimal solution to $\boldsymbol{\theta}$ is the maximum likelihood estimate, $\hat{\boldsymbol{\theta}}_{\mathrm{ML}}$, which can be obtained by minimizing the negative log-likelihood function, i.e.,

$$\hat{\boldsymbol{\theta}}_{\mathrm{ML}} \triangleq \min_{\boldsymbol{\theta}} \sum_{i=1}^{n} f_i(\boldsymbol{\theta}) \quad (2)$$

where $f_i(\boldsymbol{\theta}) \triangleq \frac{1}{2} (\mathbf{z}_i - \mathbf{h}_i(\boldsymbol{\theta}))^T R_i^{-1} (\mathbf{z}_i - \mathbf{h}_i(\boldsymbol{\theta}))$. The objective is to solve (2) using only local interactions dictated by the network topology, i.e., each agent recovers the global minimizer $\hat{\boldsymbol{\theta}}_{\mathrm{ML}}$ by only assuming access to local information $f_i(\cdot)$ and communication to one-hop neighbors.

## IV. MAIN RESULTS

Assuming $\mathbf{h}_i$'s are continuously differentiable, a first-order algorithm of the following form can be utilized to solve (2) in a centralized manner:

$$\hat{\boldsymbol{\theta}}^{k+1} = \hat{\boldsymbol{\theta}}^k - \alpha \sum_{i=1}^n \nabla f_i\left(\hat{\boldsymbol{\theta}}^k\right), \quad (3)$$

where $\nabla f_i\left(\hat{\boldsymbol{\theta}}^k\right) = \left(H_i^k\right)^T R_i^{-1}\left(\mathbf{h}_i\left(\hat{\boldsymbol{\theta}}^k\right) - \mathbf{z}_i\right)$ denotes the gradient, $\alpha > 0$ is the step size, and $H_i^k = \frac{\partial \mathbf{h}_i}{\partial \boldsymbol{\theta}}(\hat{\boldsymbol{\theta}}^k)$. Distributed implementation of the algorithm in (3) requires each agent to first reach consensus on an initial value $\hat{\boldsymbol{\theta}}^0$ and calculate the global sum $\sum_{i=1}^n \nabla f_i\left(\hat{\boldsymbol{\theta}}^k\right)$ via local interactions. Here we propose distributed optimization algorithm that utilizes a static average consensus algorithm to asymptotically reach consensus on the initial value and a DAC algorithm to reach consensus on the descent direction in finite-time.

### A. Dynamic Average Consensus

In continuous-time formulation, (3) can be written as

$$\dot{\hat{\boldsymbol{\theta}}}(t) = -\alpha \sum_{i=1}^n \nabla f_i\left(\hat{\boldsymbol{\theta}}(t)\right), \quad \hat{\boldsymbol{\theta}}(t_0) = \hat{\boldsymbol{\theta}}^0. \quad (4)$$

Let $\phi_i(t) = \nabla f_i\left(\hat{\boldsymbol{\theta}}_i(t)\right) \in \mathbb{R}^r$, where $\hat{\boldsymbol{\theta}}_i(t)$'s are the local estimates of $\hat{\boldsymbol{\theta}}_{\text{ML}}$. If all the agents have the same initial condition, and they are able to reach consensus on $\sum_{i=1}^n \phi_i(t)$, then the optimization iteration can be implemented distributedly with each local estimates tracking the "optimal" central estimate trajectory. In this subsection we propose a robust DAC algorithm that would allow the agents to reach consensus on the average $\bar{\phi}(t) = \frac{1}{n}\sum_{i=1}^n \phi_i(t)$.

Let $x_i(t) \in \mathbb{R}^r$ denote node $i$'s estimate of $\bar{\phi}(t)$. Here we propose the following DAC algorithm:

$$\dot{\mathbf{z}}(t) = \mathbf{u}(t), \quad \mathbf{z}(t_0) = \mathbf{z}_0, \quad (5a)$$
$$\mathbf{x}(t) = (\mathcal{B} \otimes I_r)\mathbf{z}(t) + \boldsymbol{\phi}(t), \quad (5b)$$

where $\boldsymbol{\phi}(t) \triangleq \left[\phi_1^T(t) \ldots \phi_n^T(t)\right]^T \in \mathbb{R}^{nr}$, $\mathbf{x}(t) \triangleq \left[x_1^T(t) \ldots x_n^T(t)\right]^T \in \mathbb{R}^{nr}$, $\mathbf{z}(t) \in \mathbb{R}^{r\ell}$ is the internal state of the entire network, and $\mathbf{u}(t)$ is the input that needs to be designed.

The following Theorem illustrate how to select the inputs $\mathbf{u}(t)$ such that the DAC-error $\tilde{\mathbf{x}}(t) \triangleq \mathbf{x}(t) - \mathbf{1}_n \otimes \bar{\phi}(t)$ converges to zero in finite time.

*Theorem 1:* For any connected undirected network, given $\sup_{t\in[t_0,\infty)}\|\dot{\boldsymbol{\phi}}(t)\|_\infty < \infty$, the robust DAC algorithm (5) guarantees that the average consensus error, $\tilde{\mathbf{x}}(t)$, converges to zero in finite time for any initial condition $\mathbf{z}_0$, if the estimator input $\mathbf{u}(t)$ is selected as

$$\mathbf{u}(t) = -\beta \operatorname{sgn}\left\{\left(\mathcal{B}^T \otimes I_r\right)\mathbf{x}(t)\right\}, \quad (6)$$

where $\beta$ is the input gain. More specifically, we have $\tilde{\mathbf{x}}(t) = 0$ for all $t \geq t^*$, where

$$t^* = \frac{\|\tilde{\mathbf{x}}(t_0)\|_2}{\sqrt{\lambda_2(\mathcal{L})}}. \quad (7)$$

*Proof:* Note $\tilde{\mathbf{x}}(t) = (\mathcal{B} \otimes I_r)\mathbf{z}(t) + (M \otimes I_r)\boldsymbol{\phi}(t)$. Thus, $\dot{\tilde{\mathbf{x}}}(t) = (\mathcal{B} \otimes I_r)\dot{\mathbf{z}}(t) + (M \otimes I_r)\dot{\boldsymbol{\phi}}(t)$. Consider a nonnegative function of the form $V = \frac{1}{2}\tilde{\mathbf{x}}^T(t)\tilde{\mathbf{x}}(t)$. Therefore,

$$\dot{V} = \tilde{\mathbf{x}}^T(t)(\mathcal{B} \otimes I_r)\mathbf{u}(t) + \tilde{\mathbf{x}}^T(t)(\mathcal{B} \otimes I_r)$$
$$\times \left(\mathcal{B}^T\left(\mathcal{B}\mathcal{B}^T\right)^+ \otimes I_r\right)\dot{\boldsymbol{\phi}}(t),$$
$$\leq \tilde{\mathbf{x}}^T(t)(\mathcal{B} \otimes I_r)\mathbf{u}(t) + \|\tilde{\mathbf{x}}^T(t)(\mathcal{B} \otimes I_r)\|_1$$
$$\times \|\left(\mathcal{B}^T\left(\mathcal{B}\mathcal{B}^T\right)^+ \otimes I_r\right)\|_\infty \|\dot{\boldsymbol{\phi}}(t)\|_\infty$$
$$\leq \tilde{\mathbf{x}}^T(t)(\mathcal{B} \otimes I_r)\mathbf{u}(t) + \|\tilde{\mathbf{x}}^T(t)(\mathcal{B} \otimes I_r)\|_1$$
$$\times \|\left(\mathcal{B}^T \otimes I_r\right)\|_\infty \|\left(\left(\mathcal{B}\mathcal{B}^T\right)^+ \otimes I_r\right)\|_\infty \|\dot{\boldsymbol{\phi}}(t)\|_\infty.$$

Substituting (6) yields

$$\dot{V} \leq -\beta\|\tilde{\mathbf{x}}^T(t)(\mathcal{B} \otimes I_r)\|_1 + \|\tilde{\mathbf{x}}^T(t)(\mathcal{B} \otimes I_r)\|_1$$
$$\times \|\left(\mathcal{B}^T \otimes I_r\right)\|_\infty \|\left(\left(\mathcal{B}\mathcal{B}^T\right)^+ \otimes I_r\right)\|_\infty \|\dot{\boldsymbol{\phi}}(t)\|_\infty.$$

Note $\|\left(\mathcal{B}^T \otimes I_r\right)\|_\infty = \|\mathcal{B}^T\|_\infty \leq 2n$ and $\|\left(\mathcal{B}\mathcal{B}^T\right)^+\|_\infty \leq \frac{\sqrt{n}}{2\lambda_2(\mathcal{L})}$, where $\lambda_2(\mathcal{L})$ denotes the algebraic connectivity. Thus, if $\beta$ is selected such that

$$\beta \geq \frac{n\sqrt{n}}{\lambda_2(\mathcal{L})}\dot{\phi}_{\max} + 1,$$

where $\dot{\phi}_{\max} \geq \sup_{t\in[t_0,\infty)}\|\dot{\boldsymbol{\phi}}(t)\|_\infty$, then we have

$$\dot{V} \leq -\|\left(\mathcal{B}^T \otimes I_r\right)\tilde{\mathbf{x}}(t)\|_1 \leq -\sqrt{\|\left(\mathcal{B}^T \otimes I_r\right)\tilde{\mathbf{x}}(t)\|_2^2}$$
$$\leq -\sqrt{\tilde{\mathbf{x}}^T(t)(\mathcal{B} \otimes I_r)(\mathcal{B}^T \otimes I_r)\tilde{\mathbf{x}}(t)}$$
$$= -\sqrt{2\tilde{\mathbf{x}}^T(t)(\mathcal{L} \otimes I_r)\tilde{\mathbf{x}}(t)} \leq -\sqrt{2}\sqrt{\lambda_2(\mathcal{L})}\sqrt{V}.$$

Thus we have $\frac{1}{2\sqrt{V}}\dot{V} \leq -\frac{1}{2}\sqrt{2\lambda_2(\mathcal{L})}$. Now based on the Comparison Lemma (Lemma 3.4 of [36]), $\sqrt{V(t)} \leq \sqrt{V(t_0)} - \frac{1}{2}\sqrt{2\lambda_2(\mathcal{L})}\,t$. Since $\dot{V}(t)$ is negative definite and $V(t)$ is positive definite, we have $\tilde{\mathbf{x}}(t) = \mathbf{0}$ for all $t \geq t^*$, where $t^* = \frac{\|\tilde{\mathbf{x}}(t_0)\|_2}{\sqrt{\lambda_2(\mathcal{L})}}$. This concludes the proof. ∎

## B. Distributed Algorithm

Based on the results obtained from the previous subsection, we propose the following distributed algorithm to solve the optimization problem given in (2):

$$\dot{\Theta}(t) = -\gamma \left( \mathcal{L} \otimes I_r \right) \Theta(t) - \hat{\alpha} \mathbf{x}(t), \quad \Theta(t_0) \quad (8a)$$
$$\dot{\mathbf{z}}(t) = -\beta \mathrm{sgn}\left\{ \left( \mathcal{B}^T \otimes I_r \right) \mathbf{x}(t) \right\}, \quad \mathbf{z}(t_0) = \mathbf{z}_0, \quad (8b)$$
$$\mathbf{x}(t) = \left( \mathcal{B} \otimes I_r \right) \mathbf{z}(t) + \phi(t), \quad (8c)$$

where $\gamma$ is the static average consensus gain, $\hat{\alpha}$ is the step size, $\Theta(t) \triangleq \left[ \hat{\boldsymbol{\theta}}_1^T(t), \ldots, \hat{\boldsymbol{\theta}}_n^T(t) \right]^T \in \mathbb{R}^{nr}$ is the concatenated vector of all local estimates of $\hat{\boldsymbol{\theta}}_{\mathrm{ML}}$. Before we present the main result of the paper, we first provide the following Lemma regarding the exponential convergence of static average consensus algorithm.

*Lemma 1:* Let the vector of the values of all the nodes $\mathbf{y}(t) \in \mathbb{R}^n$ be the solution of the following differential equation:

$$\dot{\mathbf{y}}(t) = -\mathcal{L}\mathbf{y}(t), \quad \mathbf{y}(t_0) = \mathbf{y}_0 \quad (9)$$

Then, all the nodes of the graph globally, asymptotically reach an average consensus value $\bar{y} = \frac{1}{n}\mathbf{1}_n^T \mathbf{y}_0$ with an exponential rate of $\kappa = \lambda_2(\mathcal{L})$, i.e.,

$$\| \boldsymbol{\delta}(t) \| \leq \| \boldsymbol{\delta}(t_0) \| \exp(-\kappa t), \quad (10)$$

where $\boldsymbol{\delta}(t) = \mathbf{y}(t) - \bar{y}\mathbf{1}_n$.

*Proof:* The proof follows from noticing that the solution to (9) can be written as $\mathbf{y}(t) = \exp\{-\mathcal{L}t\} \mathbf{y}_0$, and $\lim_{t \to \infty} \exp\{-\mathcal{L}t\} = \frac{1}{n}\mathbf{1}_n\mathbf{1}_n^T$. ∎

Note that the consensus protocol given in (9) is same as the gradient-descent algorithm for minimizing the *Laplacian potential* $\Psi_{\mathcal{G}}(\mathbf{y}) = \frac{1}{2}\mathbf{y}^T \mathcal{L} \mathbf{y}$ and the exponential convergence rate of the protocol can be arbitrarily increased by simply multiplying the Laplacian matrix by a positive constant $\gamma > 1$. The local implementation of (8a) for all $i \in \mathcal{I}$ can be written as

$$\dot{\hat{\boldsymbol{\theta}}}_i(t) = -\gamma \sum_{j: v_j \in \mathcal{N}_i} \left( \hat{\boldsymbol{\theta}}_i(t) - \hat{\boldsymbol{\theta}}_j(t) \right) - \hat{\alpha}\, x_i(t). \quad (11)$$

The main objective here is to make sure that the trajectories of (11) track the trajectories of (4) as close as possible to ensure that the individual local solutions to (2) asymptotically recovers the centralized performance. Now we present the main contribution of the paper which combines the results from Lemma 1 and Theorem 1.

*Theorem 2:* For any connected undirected network, the distributed algorithm in (8) converges to the centralized solution trajectories of (4) exponentially fast for all $t \geq t^*$, where $t^*$ is given in (7).

*Proof:* Based on the finite time convergence results of Theorem 1, for $t \geq t^*$, (11) can be written as

$$\dot{\hat{\boldsymbol{\theta}}}_i(t) = \gamma \sum_j (\hat{\boldsymbol{\theta}}_j(t) - \hat{\boldsymbol{\theta}}_i(t)) - \frac{\hat{\alpha}}{n} \sum_{l=1}^n \nabla f_l \left( \hat{\boldsymbol{\theta}}_l(t) \right). \quad (12)$$

Note that between the two terms present in the above equation, only the diffusion-term, i.e., $\gamma \sum_j \left( \hat{\boldsymbol{\theta}}_i(t) - \hat{\boldsymbol{\theta}}_j(t) \right)$, differer from agent to agent. Since the gradient-term is identical for all agents, based on Lemma 1, one can conclude that each $\hat{\boldsymbol{\theta}}_i(t)$ converges to $\frac{1}{n}\sum_{i=1}^n \hat{\boldsymbol{\theta}}_i(t)$ and therefore (12) converges to the trajectories of (4) exponentially fast, given $\hat{\boldsymbol{\theta}}^0 = \frac{1}{n}\sum_{i=1}^n \hat{\boldsymbol{\theta}}_i(t_0)$ and $\hat{\alpha} = n\,\alpha$. ∎

## V. NUMERICAL RESULTS

For numerical simulations, we consider the problem of distributed event localization using acoustic sensor network [2]. Each sensor consist of an array of microphones that can obtain the direction of arrival of the acoustic signal. Thus the measurement model is given as

$$z_i = \arctan\left( \frac{T_y - S_i^y}{T_x - S_i^x} \right) + w_i,$$

where $(T_x, T_y)$ denotes the unknown event location and the (locally) known two-dimensional, sensor locations are given as $(S_i^x, S_i^y)$.

Here $\boldsymbol{\theta} = [T_x, T_y]^T = [200, -300]^T$, $n = 7$, and $R_i = 10^{-2}, \forall i \in \{1, \ldots, 7\}$. The sensor locations and the network topology are given in Fig. 1(a). In Fig. 1(a), blue circles denote the sensor locations, the solid (black) lines between the sensors indicate communication links, and the true event location is marked using a red star. We conducted a $10^3$ Monte Carlo simulation to evaluate the performance of the proposed distributed algorithm. Figure 1(b) contains the mean-square estimation error (MSEE) for the centralized estimator[3] given in (4) (denoted as Cent.) and distributed estimator[4] given in (8) (denoted as Dist.). The mean-square tracking error (MSTE$(t)$) given in Fig. 1(c) is the difference between the centralized solution trajectory and the distributed solution trajectories, i.e., $\mathrm{MSTE}_i(t) = \frac{1}{10^3}\sum_{l=1}^{10^3} \left( \hat{\boldsymbol{\theta}}(t) - \hat{\boldsymbol{\theta}}_i(t) \right)^T \left( \hat{\boldsymbol{\theta}}(t) - \hat{\boldsymbol{\theta}}_i(t) \right)$. Finally

---

[2]Typical applications include gunfire detection and shooter localization, see [37], [38] for more details.

[3]MSEE$(t) = \frac{1}{10^3}\sum_{l=1}^{10^3} \left( \boldsymbol{\theta} - \hat{\boldsymbol{\theta}}(t) \right)^T \left( \boldsymbol{\theta} - \hat{\boldsymbol{\theta}}(t) \right)$

[4]MSEE$(t) = \frac{1}{10^3}\sum_{l=1}^{10^3} \frac{1}{7}\sum_{i=1}^{7} \left( \boldsymbol{\theta} - \hat{\boldsymbol{\theta}}_i(t) \right)^T \left( \boldsymbol{\theta} - \hat{\boldsymbol{\theta}}_i(t) \right)$

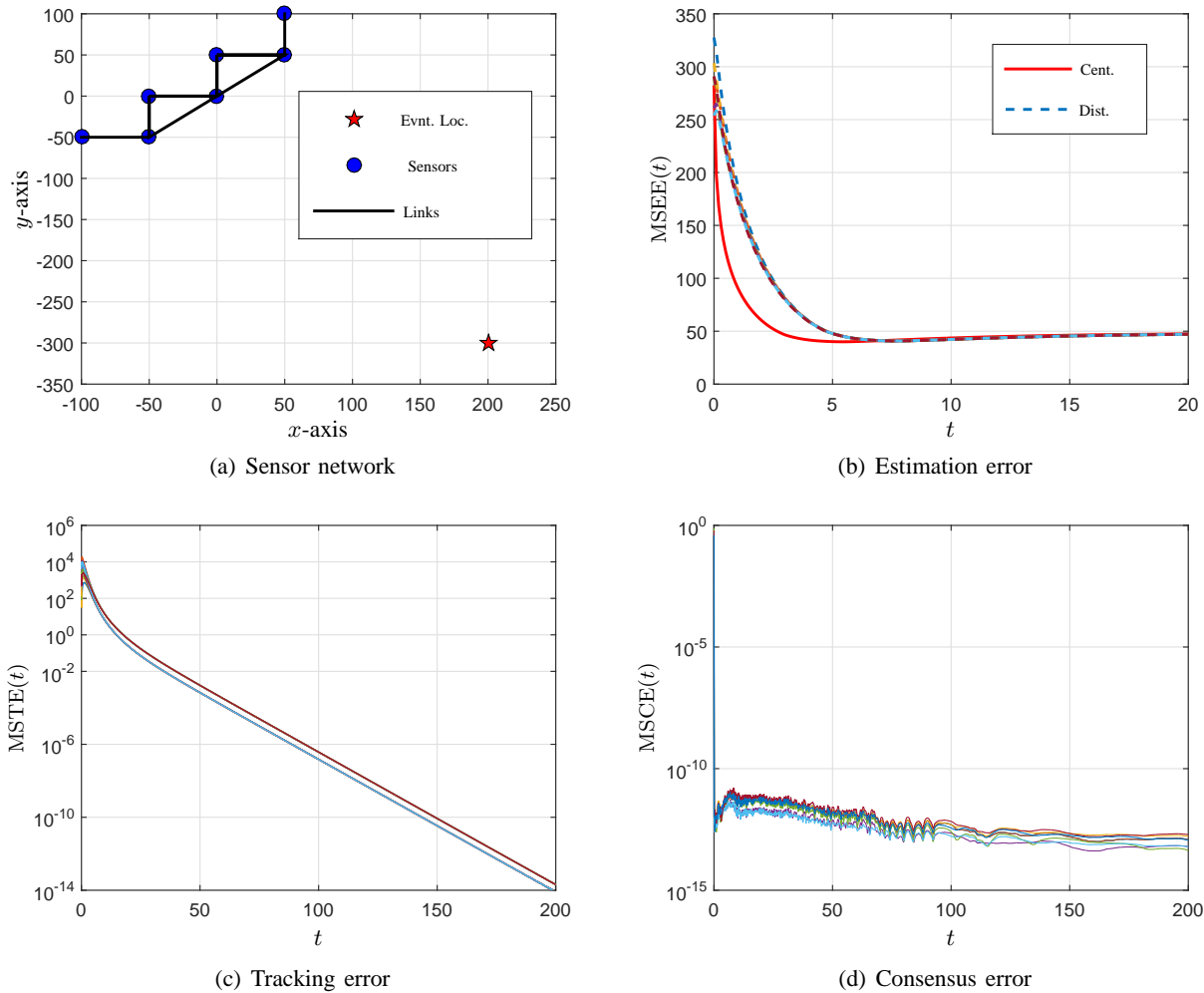

Fig. 1. Simulation scenario and mean-square errors obtained from $10^3$ Monte Carlo runs.

Fig. 1(d) contains the mean-square consensus error[5] (MSCE$(t)$) obtained from the Monte Carlo runs. Results given in Figs. 1(c) and 1(d) confirms the exponential convergence of the tracking error and the finite-time convergence of the consensus error.

## VI. CONCLUSION

This paper presents a distributed maximum likelihood scheme that utilizes a DAC algorithm to reach agreement on the gradient direction. The proposed distributed algorithm recovers the centralized performance accuracy exponentially fast. Though the current formulation focuses on first-order optimization algorithm, it can be easily applied to higher order schemes by utilizing the DAC algorithm for reaching consensus on the higher-order derivatives of the local objective function. The proposed continuous-time formulation can be extended to discrete-time scenarios after replacing the signum function with an appropriate continuous approximation such as a saturation function. Future research include extending the current approach to accelerated gradient methods and considering privacy preserving event-triggered communication schemes.

[5]MSCE$_i(t) = \frac{1}{10^3}\sum_{l=1}^{10^3} \frac{1}{7}\sum_{i=1}^{7} \left(x_i(t) - \bar{\phi}(t)\right)^T \left(x_i(t) - \bar{\phi}(t)\right)$